\documentclass[12pt, a4paper, twoside]{article}
\usepackage{amssymb}
\usepackage[leqno]{amsmath}
\usepackage{latexsym}
\usepackage{enumerate}
\usepackage{amsmath, amssymb, amscd}
\usepackage{fancyhdr}
\DeclareMathOperator{\Pic}{Pic}
\DeclareMathOperator{\Cliff}{Cliff}
\DeclareMathOperator{\rk}{rk}
\DeclareMathOperator{\Hom}{Hom}

\DeclareMathOperator{\ev}{ev}
\DeclareMathOperator{\gr}{gr}
\DeclareMathOperator{\length}{length}

\begin{document}

\title{\textbf{\Large{Slope semistability of rank 2 Lazarsfeld-Mukai bundles on K3 surfaces and ACM line bundles}}}


\date{}

\maketitle


\smallskip

\smallskip


\begin{abstract}

Previously, many people have studied a stability of vector bundles of given rank and Chern classes on algebraic varieties. Recently, we are interested in the slope stability of the rank 2 Lazarsfeld-Mukai bundle $E_{C,Z}$ on a K3 surface $X$ associated to a very ample smooth curve $C$ on $X$ and a base point free pencil $Z$ on $C$ with respect to $\mathcal{O}_X(C)$. In this paper, we will give a sufficient condition for such a Lazarsfeld-Mukai bundle $E_{C,Z}$ to be $\mathcal{O}_X(C)$-slope semistable by ACM line bundles with respect to $\mathcal{O}_X(C)$.

\end{abstract}

\section{Introduction}

The study of the Gieseker (or slope) stability and the notion of arithmetically Cohen-Macaulay (ACM for short) of vector bundles with respect to a given very ample line bundle on algebraic varieties is a very active topic in algebraic geometry. Recently, several authors have studied about the Gieseker stability of ACM bundles (that is, vector bundles with no intermediate cohomology). For example, Marta Casanellas and Robin Hartshorne ([C-H]) have proved the existence of stable Ulrich bundles of given rank and Chern classes on smooth cubic surfaces in $\mathbb{P}^3$ and investigated the structure of the moduli space of them. Emre Coskun, Rajesh S. Kulkarni and Yusuf Mustopa ([C-K-M]) have proved that all smooth quartic surfaces in  $\mathbb{P}^3$ admit a 14-dimensional family of simple Ulrich bundles of rank 2 with $c_1=3H$ and $c_2=14$, where $H$ is the very ample line bundle given by a hyperplane section. However, in general, it is difficult to investigate the stability of ACM bundles with respect to a given polarization.

In this paper, we will focus on the slope stability of a Lazarsfeld-Mukai bundle on a K3 surface with respect to the first Chern class of it. Let $X$ be a K3 surface, $C$ be a smooth curve on $X$, and $Z$ be a base point free line bundle on $C$. Then, the Lazarsfeld-Mukai bundle $E_{C,Z}$ associated to $C$ and $Z$ is defined as the dual of the kernel of the evaluation map associated to the space of the global sections of $Z$, and satisfies $h^1(E_{C,Z})=h^2(E_{C,Z})=0$. 
In particular, we can easily see that if $C$ is very ample and $|Z|$ is a pencil, that is, $E_{C,Z}$ is rank 2, then it is ACM with respect to $\mathcal{O}_X(C)$ (see Lemma 5.2). Therefore, we are interested in the $\mathcal{O}_X(C)$-slope stability of such a Lazarsfeld-Mukai bundle $E_{C,Z}$. As a previous result about the case where $E_{C,Z}$ is rank 2, Margherita Lelli Chiesa ([LC]) have proved that if $C$ is a $\lfloor\frac{g+3}{2} \rfloor$-gonal curve of genus $g$ and Clifford dimension one and $d=\deg(Z)$ satisfies $\rho(g,1,d)=2d-g-2>0$, then  $E_{C,Z}$ is slope stable with respect to $\mathcal{O}_X(C)$. However, in other cases, there is no concrete description in terms of the condition for $E_{C,Z}$ to be $\mathcal{O}_X(C)$-slope stable. In this paper, we will show that if $E_{C,Z}$ is not slope semistable with respect to $\mathcal{O}_X(C)$, the maximal destabilizing sheaf of it contains an initialized and ACM line bundle with respect to $\mathcal{O}_X(C)$ to give a sufficient condition for $E_{C,Z}$ to be $\mathcal{O}_X(C)$-slope semistable, and give some examples of it.




Our plan of this paper is as follows: In Section 2, we recall some fundamental results about linear systems on K3 surfaces. In Section 3, we recall some basic notions and results about the Clifford index of a smooth curve on a K3 surface. In Section 4, we recall the notion of the slope (semi)stability of vector bundles and ACM bundles, and prepare a proposition to prove our main result. In section 5, we recall the properties of Lazarsfeld-Mukai bundles. In Section 6, we give our main result and some examples of it.
 
$\;$

\noindent {\bf{Notations and Conventions.}} We work over the complex number field $\mathbb{C}$. A surface and a curve are smooth projective. A K3 surface is a regular surface whose canonical line bundle is trivial. 

For a curve $C$, we denote by $K_C$ the canonical line bundle of $C$. We denote by $g^r_d$ a linear system of dimension $r$ and degree $d$. For a line bundle $Z$ on a curve $C$, we denote by $|Z|$ the linear system defined by $Z$. A curve $C$ is called {\it{$k$-gonal}} if $C$ has a $g^1_k$ but no $g^1_{k-1}$. In particular, a 2-gonal curve is called {\it{hyperelliptic}}. If a curve $C$ is $k$-gonal, then $g^1_k$ is base point free and complete. Note that if $C$ is a very ample curve on a K3 surface, $C$ is not hyperelliptic (see [SD], Theorem 5.2). For a line bundle $Z$ on a curve $C$, the {\it{Clifford index}} of $Z$ is defined as follows;
$$\Cliff(Z):=\deg Z-2\dim |Z|.$$
\noindent Moreover, the {\it{Clifford index of $C$}} is defined as follows;
$$\Cliff(C):=\min\{\Cliff(Z)\;|\;h^0(Z)\geq2,\;h^1(Z)\geq2\}.$$
\noindent It is well known that one gets $0\leq \Cliff(C)\leq \lfloor\frac{g-1}{2} \rfloor$ from Brill-Noether theory (cf, [A-C-G-H, V]). Moreover, if a curve $C$ of genus $g\geq2$ is $k$-gonal, then the following inequality holds (cf. [C-M]).
$$\Cliff(C)+2\leq k\leq \Cliff(C)+3.$$
A $k$-gonal curve $C$ satisfying $k=\Cliff(C)+3$ is called an {\it{exceptional curve}}.

For two divisors $D_1$ and $D_2$ on a surface, we will write $D_1\sim D_2$ if they are linearly equivalent. For a torsion free sheaf $E$, we denote by $E^{\vee}$ the dual of $E$. A vector bundle $E$ is called {\it{simple}} if $\Hom(E,E)\cong\mathbb{C}$.

\section{Linear systems on K3 surfaces} In this section, we recall a few classical results about divisors and line bundles on K3 surfaces. 

\newtheorem{df}{Definition}[section]

\begin{df} {\rm{A divisor}} $D$ {\rm{on a surface is called}} $m$-connected {\rm{if}} $D_1.D_2\geq m$, {\rm{for each effective decomposition}} $D=D_1+D_2$.\end{df}

\noindent We can easily see that if the linear system $|L|$ defined by a line bundle $L$ on a K3 surface contains a 1-connected divisor, then $h^1(L)=0$ (for example, see [B-P-V], Corollary 12.3). On the other hand, one has the following characterization of numerical effective divisors and base point free divisors.

\newtheorem{prop}{Proposition}[section]

\begin{prop}{\rm{([SD], Proposition 2.7)}} Let $D$ be a non-zero numerical effective divisor on a K3 surface $X$. Then $D$ is not base point free if and only if there exists an elliptic curve $F$, a smooth rational curve $\Gamma$ and an integer $k\geq2$ such that $F.\Gamma=1$ and $D\sim kF+\Gamma$. \end{prop}

\begin{prop}{\rm{([SD], Proposition 2.6)}} Let $D$ be a non-zero effective divisor on a K3 surface $X$. Assume that $|D|$ has no fixed components. Then one of the following cases occurs.

\smallskip

\smallskip

{\rm{(i)}} $D^2>0$ and the general member of $|D|$ is a smooth irreducible curve of genus $\frac{1}{2}D^2+1$.

{\rm{(ii)}} $D^2=0$ and $D\sim kF$, where $k\geq1$ and $F$ is a smooth curve of genus one. In this case, $h^1(\mathcal{O}_X(D))=k-1$. \end{prop} 

\noindent Proposition 2.2 is called a strong Bertini's theorem. Moreover, since the linear system $|C|$ defined by an irreducible curve $C$ on a K3 surface with $C^2>0$ is base point free ([SD], Theorem 3.1), one also has the following assertion as a corollary of Proposition 2.2.

\begin{prop}{\rm{([SD], Corollary 3.2)}} Let $L$ be a line bundle on a K3 surface. Then $|L|$ has no base points outside its fixed components.\end{prop}

\noindent In particular, a very ample line bundle is not hyperelliptic. Hence, by the characterization of hyperellitic linear systems (cf. [M-M], and [SD], Theorem 5.2), we have the following assertion.

\begin{prop} Let $L$ be a numerical effective line bundle with $L^2\geq4$ on a K3 surface $X$. Then $L$ is very ample if and only if the following conditions are satisfied.

\smallskip

\smallskip

{\rm{(i)}} There is no irreducible curve $E$ such that $E^2=0$ and $E.L=1$ or 2.

{\rm{(ii)}} There is no irreducible curve $E$ such that $E^2=2$ and $L\cong\mathcal{O}_X(2E)$.

{\rm{(iii)}} There is no irreducible curve $E$ such that $E^2=-2$ and $E.L=0$. \end{prop}

\section{The Clifford index of smooth curves on K3 surfaces} In this section, we recall some results about the Clifford index of curves on K3 surfaces. First of all, we prepare some notations to explain them.

\begin{df} {\rm{Let}} $X$ {\rm{be a K3 surface, and let}} $L$ {\rm{be a base point free and big line bundle on}} $X$. {\rm{Then let}}
$$\mathcal{A}(L):=\{D\in\Pic(X)\;|\;h^0(D)\geq2,\;h^0(L-D)\geq2\}.$$
\noindent {\rm{If}} $\mathcal{A}(L)$ {\rm{is not empty, let}}
$$\mu(L):=\min\{D.(L-D)-2\;|\;D\in\mathcal{A}(L)\},$$
\noindent {\rm{and set}}
$$\mathcal{A}^0(L):=\{D\in\mathcal{A}(L)\;|\;D.(L-D)=\mu(L)+2\}.$$ 
\end{df}

\noindent For a line bundle $L$ as in Definition 3.1 and divisors belonging to $\mathcal{A}^0(L)$, Johnsen and Knutsen showed the following result (cf. [J-K] Proposition 2.6). 

\begin{prop} If $\mathcal{A}(L)\neq\emptyset$, then $\mu(L)\geq0$, and all divisors $D\in \mathcal{A}^0(L)$ satisfy the following conditions.

\smallskip

\smallskip

\noindent {\rm{(i)}} The base divisor $\Delta$ of $D$ satisfies $L.\Delta=0$.

\noindent {\rm{(ii)}} $h^1(D)=0.$

\end{prop}

\noindent In Proposition 3.1, if $\mathcal{A}(L)\neq\emptyset$, then we can find $D\in\mathcal{A}^0(L)$ such that either $|D|$ or $|L-D|$ is base point free and its general member is an irreducible smooth curve (cf. [J-K] Proposition 2.7).

\newtheorem{thm}{Theorem}[section]

\begin{thm} {\bf{[G-L]}} Let $X$ be a K3 surface, and let $L$ be a base point free and big line bundle of sectional genus $g$ on $X$. Then the Clifford index of the smooth curves of $|L|$ is constant, and for any smooth curve $C\in |L|$, if $\Cliff(C)< \lfloor\frac{g-1}{2} \rfloor$, then there exists a divisor $D$ on $X$ such that $\Cliff(\mathcal{O}_C(D))=\Cliff(C)$.\end{thm} 

\noindent By Theorem 3.1, the Clifford index of $L$ as in Theorem 3.1 can be defined by the Clifford index of the smooth curves of $|L|$, and it is denoted by $\Cliff(L)$. By the proof of Theorem 3.1, one can choose the divisor $D$ as in Theorem 3.1 so that it is smooth and belongs to $\mathcal{A}(L)$. Knutsen proved the following result (cf. [Kn] Lemma 8.3).

\begin{thm} Let $L$ be as in Theorem 3.1 and assume that $\Cliff(L)=c$. If $c<\lfloor\frac{g-1}{2} \rfloor$, then there exists a smooth curve $D$ on $X$ satisfying $0\leq D^2\leq c+2,\;2D^2\leq D.L$ (either of the latter two inequalities being an equality if and only if $L\sim 2D$) and 
$$c=\Cliff(\mathcal{O}_C(D))=D.L-D^2-2,$$
\noindent for any smooth curve $C\in|L|$.\end{thm}

\noindent By  Proposition 3.1, Theorem 3.1, and Theorem 3.2, for a base point free and big line bundle $L$ on a K3 surface of sectional genus $g$, we have $$\Cliff(L)=\min\{\lfloor\frac{g-1}{2} \rfloor,\mu(L)\}.$$

\section{Slope (semi)stability of vector bundles and ACM line bundles on K3 surfaces} In this section, we recall the notion of ACM bundles and the slope (semi)stability of vector bundles with respect to a given polarization on a K3 surface. Let $X$ be a surface, and let $H$ be a very ample line bundle on $X$ which provides a closed embedding in a projective space of higher dimension. Then, we denote the line bundle $H^{\otimes l}$ by $\mathcal{O}_X(l)$. For a vector bundle $E$ on a surface $X$, we will write $E(l)=E\otimes\mathcal{O}_X(l)$.

\begin{df} {\rm{A vector bundle}} $E$ {\rm{on a surface}} $X$ {\rm{is called}} initialized {\rm{if}} $H^0(E)\neq0$ {\rm{but}} $H^0(E(-1))=0$.\end{df}

\begin{df} {\rm{A vector bundle}} $E$ {\rm{on a surface}} $X$ {\rm{is called an}} Arithmetically Cohen-Macaulay {\rm{(}}ACM for short{\rm{)}} {\rm{if}} $H^1(E(l))=0$, {\rm{for all}} $l\in\mathbb{Z}$.\end{df}

\noindent In particular, for a line bundle on a K3 surface, we have the following assertion.

\begin{prop}{\rm{([W], Lemma 3.1)}}. Assume that $X$ is a K3 surface and let $L$ be a line bundle on $X$ with $|L|\neq\emptyset$. Moreover, let $m\in\mathbb{N}$. Then if $H.L\leq mH^2-1$ and, for any $k\in\mathbb{Z}$ with $0\leq k\leq m$, $h^1(L(-k))=0$, then $L$ is an ACM line bundle with respect to $H$. \end{prop}

\noindent In Proposition 4.1, we can easily see that if $L\in\mathcal{A}(H)$, $L$ is ACM precisely when $h^1(L)=h^1(H\otimes L^{\vee})=0$, since $|L|\neq\emptyset$ and $|H\otimes L^{\vee}|\neq\emptyset$. Note that such a line bundle $L$ is initialized with respect to $H$, since $h^0(L\otimes H^{\vee})=0$. 

Next, we recall the definition and some facts about the slope stability of vector bundles (cf. [HL] and [Sh]). 

\begin{df} {\rm{Let}} $X$ {\rm{and}} $H$ {\rm{be as above, and let}} $E$ {\rm{be a torsion free sheaf on}} $X$ {\rm{of rank}} $r$. {\rm{Then the}} $H${\rm{-slope of}} $E$ {\rm{is defined as follows}};
$$\mu_H(E)=\frac{c_1(E).H}{r}.$$
\noindent $E$ {\rm{is called}} $\mu_H$-semistable {\rm{(resp.}} $\mu_H$-stable{\rm{)}} {\rm{if for any subsheaf}} $0\neq F\subset E$ {\rm{with}} $\rk F<\rk E$, {\rm{we have}} $\mu_H(F)\leq\mu_H(E)$ {\rm{(resp}}. $\mu_H(F)<\mu_H(E)${\rm{)}}.\end{df}

\noindent It is well known that for a vector bundle $E$ on $X$ and a given polarization $H$, there is a unique filtration called the {\it{Harder-Narasimhan}} ({\it{HN for short}}) {\it{filtration}} 
$$0=E_0\subset E_1\subset\cdots\subset E_n=E$$
\noindent such that $E_i$ is locally free and $E_i/E_{i-1}$ ($1\leq i\leq n$) are torsion free and $\mu_H$-semistable sheaves with $\mu_H(E_{i+1}/E_i)<\mu_H(E_i/E_{i-1})$, for $1\leq i\leq n-1$. Moreover, such a filtration satisfies the following inequality
$$\mu_H(E_1)>\mu_H(E_2)>\cdots>\mu_H(E).$$
\noindent Obviously, if $E$ is not $\mu_H$-semistable, then $n\geq2$. The sheaf $E_1$ is called the {\it{maximal destabilizing sheaf}} of $E$. Moreover, if a vector bundle $E$ is $\mu_H$-semistable, there exists a filtration called a {\it{Jordan-H\"{o}lder}} ({\it{JH for short}}) {\it{filtration}}
$$0=JH_0(E)\subset JH_1(E)\subset\cdots\subset JH_m(E)=E$$
\noindent such that $\gr_i(E):=JH_i(E)/JH_{i-1}(E)$ is a torsion free and $\mu_H$-stable sheaf whose slope is equal to $\mu_H(E)$ for $1\leq i\leq m$.

\section{Structures of Lazarsfeld-Mukai bundles} In this section, we recall the definition and some properties of the Lazarsfeld-Mukai bundle associated to a smooth curve on a K3 surface and a base point free line bundle on it, and prepare some lemmas to explain the main result in the next section.
$\;$

Let $X$ be a K3 surface, let $C$ be a smooth curve of genus $g\geq2$ on $X$, and let $Z$ be a base point free divisor on $C$. Then, for the evaluation map $$\ev_{Z,X}:H^0(\mathcal{O}_C(Z))\otimes\mathcal{O}_X\longrightarrow \mathcal{O}_C(Z),$$
\noindent we set $F_{C,Z}:=\ker(\ev_{Z,X})$ and $E_{C,Z}:=F_{C,Z}^{\vee}$. Since $Z$ is base point free, $F_{C,Z}$ is a locally free sheaf which fits into the following exact sequence;
$$0\longrightarrow F_{C,Z}\longrightarrow H^0(\mathcal{O}_C(Z))\otimes\mathcal{O}_X\xrightarrow{\ev_{Z,X}}\mathcal{O}_C(Z)\longrightarrow0.$$
\noindent Taking the dual of it, we get
$$0\longrightarrow H^0(\mathcal{O}_C(Z))^{\vee}\otimes\mathcal{O}_X\longrightarrow E_{C,Z}\longrightarrow K_C\otimes\mathcal{O}_C(-Z)\longrightarrow0.$$
\noindent The vector bundle $E_{C,Z}$ defined as above is called a {\it{Lazarsfeld-Mukai bundle}}, and has the following properties (for example, see [L1], [L2], [P]). 

\begin{prop} If $E_{C,Z}$ is the Lazarsfeld-Mukai bundle associated to a smooth curve $C$ of genus $g$ on a K3 surface $X$ and a base point free divisor $Z$ on $C$ of degree $d$ such that $\dim|Z|=r$, then we get the following assertion. 

$\;$

{\rm{(a)}} $c_1(E_{C,Z})=\mathcal{O}_X(C).$

\smallskip

{\rm{(b)}} $c_2(E_{C,Z})=d.$

\smallskip

{\rm{(c)}} $h^2(E_{C,Z})=h^1(E_{C,Z})=0.$

\smallskip

{\rm{(d)}} $E_{C,Z}$ is globally generated off the base points of $K_C\otimes\mathcal{O}_C(-Z)$.

\smallskip

{\rm{(e)}} If $\rho(g,r,d)=g-(r+1)(g-d+r)<0$, then $E_{C,Z}$ is non-simple.

\end{prop}

\noindent By Proposition 5.1 (d), if $|K_C\otimes\mathcal{O}_C(-Z)|\neq\emptyset$, $E_{C,Z}$ is globally generated off a finite set. Hence, we have the following assertion.

\begin{prop} {\rm{(proof of [LC], Lemma 3.2)}} Let $X$ and $E_{C,Z}$ be as above, and $Q$ be a torsion free sheaf of rank 1 on $X$. Assume that $|K_C\otimes\mathcal{O}_C(-Z)|\neq\emptyset$. If there exists a surjective morphism $\varphi:E_{C,Z}\rightarrow Q$, then $Q^{\vee\vee}$ is base point free and not trivial.\end{prop}

\noindent If $E_{C,Z}$ is rank two, that is, $|Z|$ is a pencil, one can get the following characterization and a more detailed property of it.

\newtheorem{lem}{Lemma}[section]

\begin{lem} {\bf{([D-M], Lemma 4.4, and [C-P], Lemma 2.1)}} If $E_{C,Z}$ is rank 2 and non-simple, then there exist a 0-dimensional subscheme $Z^{'}$ in $X$ which is a locally complete intersection and two line bundles $M$ and $N$ on $X$ with $h^0(N),\;h^0(M)\geq2$ such that $N$ is base point free, and $E_{C,Z}$ fits into the following exact sequence;
$$0\longrightarrow M\longrightarrow E_{C,Z}\longrightarrow N\otimes\mathcal{I}_{Z^{'}}\longrightarrow 0.$$
\noindent Moreover, if $h^0(M\otimes N^{\vee})=0$, then $Z{'}=\emptyset$ and $E_{C,Z}$ splits into the direct sum of $M$ and $N$.
\end{lem}

\noindent The exact sequence as in Lemma 5.1 is called {\it{Donagi-Morrison's extension}}, and it is uniquely determined (see [A-F], Lemma 3.2). Hence, $E_{C,Z}$ splits into a direct sum of two line bundles if and only if the Donagi-Morrison's extension associated to $E_{C,Z}$ splits. 

\begin{lem} Assume that $E_{C,Z}$ is rank 2 and $C$ is very ample as a divisor on $X$. Then, $E_{C,Z}$ is ACM and initialized with respect to $\mathcal{O}_X(C)$. \end{lem}

\smallskip

\smallskip

\noindent {\it{Proof.}} Let $H:=\mathcal{O}_X(C)$. First of all, we show that $E_{C,Z}$ is ACM with respect to $H$. Since $\det(E_{C,Z})=H$ and $E_{C,Z}$ is rank 2, for any $l\in\mathbb{Z}$, we have $E_{C,Z}(l)\cong E_{C,Z}^{\vee}(l+1)$. Hence, we have
$$h^1(E_{C,Z}(l))=h^1(E_{C,Z}^{\vee}(l+1))=h^1(E_{C,Z}(-l-1)).$$
\noindent Therefore, by Proposition 5.1, it is sufficient to show that, for any $l\geq1$, $h^1(E_{C,Z}(l))=0$.

Assume that $l\geq1$. Since $h^1(\mathcal{O}_X(l))=0$, by the exact sequence
$$0\longrightarrow H^0(\mathcal{O}_C(Z))^{\vee}\otimes\mathcal{O}_X(l)\longrightarrow E_{C,Z}(l)\longrightarrow K_C^{\otimes{l+1}}\otimes\mathcal{O}_C(-Z)\longrightarrow0,$$
\noindent we show that $h^1(K_C^{\otimes{l+1}}\otimes\mathcal{O}_C(-Z))=0$. Since $|Z|$ is a pencil, we have
$$h^1(\mathcal{O}_C(Z))=\frac{1}{2}C^2+2-\deg Z,$$
\noindent and hence, we have
$$\deg Z\leq \frac{1}{2}C^2+2.$$
\noindent Since $C$ is very ample, we have $C^2\geq4$. If $C^2\geq6$, since $\deg(\mathcal{O}_C(Z)\otimes K_C^{\vee\otimes l})<0$, we have
$$h^1(K_C^{\otimes{l+1}}\otimes\mathcal{O}_C(-Z))=h^0(\mathcal{O}_C(Z)\otimes K_C^{\vee\otimes l})=0.$$
\noindent Assume that $C^2=4$. Since $|Z|$ is a pencil, $Z$ is not linearly equivalent to $K_C$. Hence, by the same reason as above, we have $h^1(K_C^{\otimes{l+1}}\otimes\mathcal{O}_C(-Z))=0$. Therefore, we have $h^1(E_{C,Z}(l))=0$, and hence, $E_{C,Z}$ is ACM. By the way of the construction of $E_{C,Z}$, we have $h^0(E_{C,Z})\neq0$. Moreover, by the exact sequence
$$0\longrightarrow H^0(\mathcal{O}_C(Z))^{\vee}\otimes\mathcal{O}_X(-1)\longrightarrow E_{C,Z}(-1)\longrightarrow \mathcal{O}_C(-Z)\longrightarrow0,$$
\noindent we have $h^0(E_{C,Z}(-1))=0$. Hence, $E_{C,Z}$ is initialized. $\hfill\square$


\section{Slope semistability of Lazarsfeld-Mukai bundles of rank 2} In this section, we state our main result and give a proof of it. Moreover, we give some examples of it.

\begin{thm} Let $X$ be a K3 surface, let $C$ be an ample curve on $X$, and let $Z$ be a base point free divisor on $C$ such that $|Z|$ is a pencil on $C$ and $|K_C\otimes\mathcal{O}_C(-Z)|\neq\emptyset$. We set $H=\mathcal{O}_X(C)$. Assume that $E_{C,Z}$ associated to $C$ and $Z$ is not $\mu_H$-semistable. Then $E_{C,Z}$ contains an initialized and ACM line bundle $L$ with respect to $H$ such that $L^2\geq2$. Moreover, if $|Z|$ is a gonality pencil on $C$, then the maximal destabilizing sheaf of $E_{C,Z}$ is ACM with respect to $H$.\end{thm}

\smallskip

\smallskip

\noindent {\it{Proof.}} Assume that $E_{C,Z}$ is not $\mu_H$-semistable. Let $L_1$ be the maximal destabilizing sheaf of it. Since $E_{C,Z}/L_1$ is a torsion free sheaf, there exists a line bundle $L_2$ and a 0-dimensional subscheme $Z^{'}$ in $X$ such that $E_{C,Z}/L_1\cong L_2\otimes\mathcal{I}_{Z^{'}}$, where $\mathcal{I}_{Z^{'}}$ is the ideal sheaf of $Z^{'}$ in $X$. Since $|K_C\otimes\mathcal{O}_C(-Z)|\neq\emptyset$, $E_{C,Z}$ is globally generated off the base points of $|K_C\otimes\mathcal{O}_C(-Z)|$. Hence, by Proposition 5.2, $|L_2|$ is base point free and not trivial. We have $L_2^2\geq0$. Since $L_1.H\geq\frac{H^2}{2}+1$, we have $L_2.H\leq\frac{H^2}{2}-1$ and hence, we have
$$L_1^2= H^2-2L_2.H+L_2^2\geq H^2-2(\frac{H^2}{2}-1)=2.$$
\noindent Therefore, $h^0(L_1)\geq2$. Since, by Lemma 5.2, $E_{C,Z}$ is ACM with respect to $H$ and $L_2=H\otimes L_1^{\vee}$, by the exact sequence
$$0\longrightarrow L_1(-1)\longrightarrow E_{C,Z}(-1)\longrightarrow L_2(-1)\otimes\mathcal{I}_{Z^{'}}\longrightarrow 0,$$
\noindent we have $h^1(L_2)=0$. 

Assume that $L_1$ is nef. Then, we show that $L_1$ is initialized and ACM with respect to $H$. It is sufficient to show that $h^1(L_1)=0$, since $L_1\in\mathcal{A}(H)$ and $h^1(L_2)=0$, by Proposition 4.1. If $|L_1|$ is base point free, then, by the Bertini's theorem, we have the assertion. Otherwise, by Proposition 2.1, there exists an elliptic curve $F$ and a $(-2)$-curve $\Gamma$ such that $F.\Gamma=1$ and $L_1\cong\mathcal{O}_X(kF+\Gamma)\;(k\geq2)$, and hence, we have the assertion.

Assume that $L_1$ is not nef. Let $\Delta$ be the fixed component of $|L_1|$ and $D\in|L_1(-\Delta)|$. We note that since $L_1^2\geq2$, we have $D\neq0$. Here, we take a $(-2)$-curve $\Gamma_1\subset\Delta$ such that $\Gamma_1.L_1<0$. Since $H$ is ample, we have $L_2.\Gamma_1\geq2$. Since $L_2$ is base point free, $L_2(\Gamma_1)$ is also base point free and big. Let $r\geq2$, and assume that $L_2(\sum_{1\leq i\leq r-1}\Gamma_i)$ is base point free and big, and $L_1(-\sum_{1\leq i\leq r-1}\Gamma_i)$ is not nef. Then there exists a $(-2)$-curve $\Gamma_r\subset \Delta-\sum_{1\leq i\leq r-1}\Gamma_i$ such that 
$$L_1(-\sum_{1\leq i\leq r-1}\Gamma_i).\Gamma_r<0.$$ 
\noindent Since $H$ is ample, we have 
$$L_2(\sum_{1\leq i\leq r-1}\Gamma_i).\Gamma_r\geq2.$$
\noindent Since $L_2(\sum_{1\leq i\leq r-1}\Gamma_i)$ is base point free, $L_2(\sum_{1\leq i\leq r}\Gamma_i)$ is also base point free and big. Therefore, by induction, there exist a finite number of $(-2)$-curves $\Gamma_1,\cdots,\Gamma_n\subset\Delta$ such that $L_2(\sum_{1\leq i\leq n}\Gamma_i)$ is base point free and big, and $L_1(-\sum_{1\leq i\leq n}\Gamma_i)$ is nef. Here, we set $$\tilde{L_1}=L_1(-\sum_{1\leq i\leq n}\Gamma_i)\text{ and }\tilde{L_2}=L_2(\sum_{1\leq i\leq n}\Gamma_i).$$

Assume that $\tilde{L_1}^2\geq2$. Then we show that $\tilde{L_1}$ is ACM and initialized with respect to $H$. It is sufficient to show that $h^1(\tilde{L_1})=h^1(\tilde{L_2})=0$, since $\tilde{L_1}\in\mathcal{A}(H)$, by Proposition 4.1. By the Bertini's theorem, the latter equation is trivial. Moreover, the first equation also holds by the same reason as above. Hence, in this case, we have the assertion. 

Assume that $\tilde{L_1}^2=0$. Since $\tilde{L_1}$ is nef, by Proposition 2.1, we can easily see that it is base point free. Hence, by Proposition 2.2, there exists an elliptic curve $F$ such that $\tilde{L_1}\cong\mathcal{O}_X(kF)\;(k\geq1)$. If $k=1$, we have $h^1(\tilde{L_1})=0$, and hence, we have
$$\chi(\tilde{L_1})=h^0(\tilde{L_1})=h^0(L_1)\geq\chi(L_1).$$
\noindent However, this contradicts to the assumption that $\tilde{L_1}^2=0$. Hence, we have $k\geq2$. Here, we note that for any $(-2)$-curve $\Gamma\subset \Delta$, we have $F.\Gamma\leq1$. In fact, if there exists a $(-2)$-curve $\Gamma_0\subset\Delta$ such that $F.\Gamma_0\geq2$, we have the contradiction
$$h^0(\tilde{L_1})=k+1<2k+1\leq \chi(\tilde{L_1}(\Gamma_0))=h^0(\tilde{L_1}(\Gamma_0)).$$
\noindent If, for any $(-2)$-curve $\Gamma\subset \Delta$, $\Gamma.F=0$, we have the contradiction $L_1^2=\Delta^2<0$. Hence, we take a $(-2)$-curve $\Gamma\subset\Delta$ such that $\Gamma.F=1$. Since $\tilde{L_1}(\Gamma)$ is nef and, for any $0\neq\Delta^{'}\subset \Delta-\Gamma$, $\tilde{L_1}(\Delta^{'}+\Gamma)$ is not nef by Proposition 2.1. Therefore, by the same argument as above, $\tilde{L_2}(-\Gamma)$ is base point free and big. Therefore, we have $h^1(\tilde{L_1}(\Gamma))=h^1(\tilde{L_2}(-\Gamma))=0$. Since $\tilde{L_1}(\Gamma)\in\mathcal{A}(H)$, by Proposition 4.1, $\tilde{L_1}(\Gamma)$ is ACM and initialized with respect to $H$, and $\tilde{L_1}(\Gamma)^2\geq2$. 

Assume that $|Z|$ is a gonality pencil on $C$, and let $d=\deg Z$ and $c=\Cliff(C)$. Since $d=c_2(E_{C,Z})$, by the exact sequence
$$0\longrightarrow L_1\longrightarrow E_{C,Z}\longrightarrow L_2\otimes\mathcal{I}_{Z^{'}}\longrightarrow 0,$$
\noindent we have $L_1.L_2\leq d$. By the same reason as above, it is sufficient to show that $h^1(L_1)=0$. Since $c+2\leq d\leq c+3$, we have
$$c\leq\mu(H)\leq L_1.L_2-2\leq c+1,$$
\noindent and hence, we have $L_1.L_2-2=\mu(H)$ or $\mu(H)+1$. If $L_1.L_2-2=\mu(H)$, since $L_1\in\mathcal{A}^0(H)$, by Proposition 3.1, the assertion holds. We consider the case where $L_1.L_2-2=\mu(H)+1$. Assume that $h^1(L_1)\neq0$. Since $L_1^2\geq2$, by the Ramanujam's theorem, there exist non-zero effective divisors $D_1$ and $D_2$ such that $D_1^2\geq0$, $D_2^2\leq0$, $D_1.D_2\leq0$, and $L_1\cong\mathcal{O}_X(D_1+D_2)$. Since $D_1\in\mathcal{A}(H)$, we have
$$D_1.L_2(D_2)=L_1.L_2-D_2.L_2+D_1.D_2\geq\mu(H)+2,$$
\noindent and hence, we have $D_1.D_2\geq D_2.L_2-1$. Since $H$ is ample, we have
$$H.D_2=D_2L_1+D_2.L_2\geq1.$$
\noindent Since $D_2.L_1=D_2^2+D_1.D_2\leq0$ and $D_2.L_2-1\leq0$, we have
$$D_1^2=D_2^2=D_1.D_2=D_2.L_2-1=0.$$
\noindent However, this contradicts to the assumption that $L_1^2\geq2$. Hence, we have $h^1(L_1)=0$. Therefore, we have the assertion. $\hfill\square$

\smallskip

\smallskip

\noindent By Theorem 6.1, we have a sufficient condition for $E_{C,Z}$ to be $\mu_H$-semistable. Here, we give some examples of Lazarsfeld-Mukai bundles which are $\mu_H$-semistable.  

\smallskip

\smallskip

\noindent \begin{prop} Let $X\subset\mathbb{P}^3$ be a smooth quartic, $C\subset X$ be a smooth hyperplane section of $X$, and let $Z$ be a divisor such that $|Z|$ is a gonality pencil on $C$ {\rm{(}}i.e., $\deg Z=3${\rm{)}}. We set $H=\mathcal{O}_X(C)$. Then, the Lazarsfeld-Mukai bundle with respect to $C$ and $Z$ is $\mu_H$-semistable. \end{prop}

\smallskip

\smallskip

\noindent {\it{Proof}}. Assume that $E_{C,Z}$ is not $\mu_H$-semistable. Let $L_1$ be the maximal destabilizing sheaf of it. Since $E_{C,Z}/L_1$ is a torsion free sheaf, $E_{C,Z}$ fits the following exact sequence;
$$0\longrightarrow L_1\longrightarrow E_{C,Z}\longrightarrow L_2\otimes\mathcal{I}_{Z^{'}}\longrightarrow 0,$$
\noindent where $L_2$ is a line bundle and $Z^{'}$ is a 0-dimensional subscheme in $X$. Since $\deg(K_C\otimes\mathcal{O}_C(-Z))=1$, by the Riemann-Roch theorem, we have $$|K_C\otimes\mathcal{O}_C(-Z)|\neq\emptyset.$$
\noindent Therefore, by the proof of Theorem 6.1, we have $L_2^2\geq0$. Since $L_1$ is the maximal destabilizing sheaf, we have 
$$L_2.H\leq\frac{H^2}{2}-1=1.$$
\noindent However, this contradicts to the ampleness of $H$. Therefore, $E_{C,Z}$ is $\mu_H$-semistable.$\hfill\square$

\smallskip

\smallskip

\noindent In particular, $E_{C,Z}$ as in Proposition 6.1 is $\mu_H$-stable. In fact, since $E_{C,Z}$ is $\mu_H$-semistable, if it is not $\mu_H$-stable, we have a JH-filtration
$$0\subset JH_1(E_{C,Z})\subset E_{C,Z}.$$
\noindent Since $JH_1(E_{C,Z})$ and $E_{C,Z}/JH_1(E_{C,Z})$ are torsion free, there exist line bundles $L,\;L^{'}$ and 0-dimensional subschemes $Z^{'},\;Z^{''}$ in $X$ such that 
$$JH_1(E_{C,Z})\cong L\otimes\mathcal{I}_{Z^{'}}\text{ and }E_{C,Z}/JH_1(E_{C,Z})\cong L^{'}\otimes\mathcal{I}_{Z^{''}}.$$
\noindent Since $E_{C,Z}$ is globally generated off the base point of $|K_C\otimes\mathcal{O}_C(-Z)|$, $L^{'}$ is base point free and not trivial by Proposition 5.2. Since $H$ is very ample, by Proposition 2.4 and the Hodge index theorem, we have $L^{'}.H\geq3$. However, since 
$$\mu_H(JH_1(E_{C,Z}))=\mu_H(E_{C,Z}),$$
\noindent we have the contradiction
$$H.L^{'}=H^2-H.L=H.L=2.$$
\noindent Hence, $E_{C,Z}$ is $\mu_H$-stable.

\begin{prop} Let $\pi:X\rightarrow\mathbb{P}^2$ be a double covering branched along a smooth sextic. Assume that $\pi^{\ast}\mathcal{O}_{\mathbb{P}^2}(1)$ is ample. Let $H=\pi^{\ast}\mathcal{O}_{\mathbb{P}^2}(3)$, let $C\in |H|$ be a smooth curve and let $Z$ be a base point free divisor on $C$ such that $|Z|$ is a pencil and $|K_C\otimes\mathcal{O}_C(-Z)|\neq\emptyset$. If $E_{C,Z}$ associated to $C$ and $Z$ is not $\mu_H$-semistable, then there exists a 0-dimensional subscheme $Z^{'}$ with $0\leq \length Z^{'}\leq3$ such that $E_{C,Z}$ fits the following exact sequence;
$$0\longrightarrow \pi^{\ast}\mathcal{O}_{\mathbb{P}^2}(2)\longrightarrow E_{C,Z}\longrightarrow \pi^{\ast}\mathcal{O}_{\mathbb{P}^2}(1)\otimes\mathcal{I}_{Z^{'}}\longrightarrow 0.$$\end{prop}

\smallskip

\smallskip

\noindent {\it{Proof}}. Since $|Z|$ is a pencil, we have
$$h^0(K_C\otimes\mathcal{O}_C(-Z))=h^1(\mathcal{O}_C(Z))=11-\deg Z.$$
\noindent Note that since $|K_C\otimes\mathcal{O}_C(-Z)|\neq\emptyset$, we have $\deg Z\leq 10$. Assume that $E_{C,Z}$ is not $\mu_H$-semistable. Let $L_1$ be the maximal destabilizing sheaf of $E_{C,Z}$, and write $E_{C,Z}/L_1\cong L_2\otimes\mathcal{I}_{Z^{'}}$, where $L_2$ is a line bundle and $Z^{'}$ is a 0-dimensional subscheme in $X$. Since $|K_C\otimes\mathcal{O}_C(-Z)|\neq\emptyset$, by the same reason as in the proof of Theorem 6.1, $L_2$ is base point free and not trivial, and hence, $L_2^2\geq0$. Assume that $L_2^2=0$. Then there exists an elliptic curve $F$ such that $L_2\cong\mathcal{O}_X(kF)\;(k\geq1)$. If we let $H^{'}=\pi^{\ast}\mathcal{O}_{\mathbb{P}^2}(1)$, we have $H^{'}.F\geq3$, otherwise, since $H^{'2}=2$ and $H^{'}.F=2$, we have the contradiction $(H^{'}-F)^2=-2$ and $H^{'}.(H^{'}-F)=0$, by the ampleness of $H^{'}$. Therefore, we have $H.L_2\geq9$, and hence, we have the contradiction
$$H.L_1\leq\frac{H^2}{2}=9.$$
Hence, we have $L_2^2>0$. Assume that $L_2^2\geq4$. By the Hodge index theorem, we have
$$(H.L_2)^2\geq72>64.$$
\noindent This means that $H.L_1\leq 9$. This contradicts to the assumption that $L_1$ is the maximal destabilizing sheaf. Hence, we have $L_2^2=2$. By the Hodge index theorem and the assumption for $L_1$, we have
$$6\leq H.L_2\leq\frac{H^2}{2}-1=8.$$
\noindent Since $3|H.L_2$, we have $H.L_2=6$, and hence, by easy computation, we have $L_2\cong H^{'}$ and $L_1\cong H^{'\otimes2}=\pi^{\ast}\mathcal{O}_{\mathbb{P}^2}(2)$. Obviously, this means that $h^1(L_1)=h^1(L_2)=0$, and hence, $L_1$ is ACM and initialized with respect to $H$. 

On the other hand, by the exact sequence
$$0\longrightarrow H^0(\mathcal{O}_C(Z))^{\vee}\otimes\mathcal{O}_X\longrightarrow E_{C,Z}\longrightarrow K_C\otimes\mathcal{O}_C(-Z)\longrightarrow 0,$$
\noindent we have $h^0(E_{C,Z})=13-\deg Z$. Since $L_1$ is a subsheaf of $E_{C,Z}$, we have
$$h^0(E_{C,Z})\geq h^0(L_1)=\frac{L_1^2}{2}+2=6,$$
\noindent and hence, we have $\deg Z\leq7$. Since $L_1.L_2=4$ and
$$\deg Z=L_1.L_2+\length Z^{'},$$
\noindent we have $0\leq\length Z^{'}\leq 3$. Therefore, we have the assertion. $\hfill\square$

\smallskip

\smallskip

\noindent The linear system $|H|$ as in Proposition 6.2 is known as a counterexample to the conjecture of Harris and Mumford that the gonality of K3 sections should be constant, and is called Donagi-Morrison's example (cf. [D-M], 2.2). We can easily see that, if we let $d$ be the gonality of $C\in|H|$, we have $d=4$ or 5 and hence, $\rho(10,1,d)<0$. Hence, in Proposition 6.2, if we assume that $|Z|$ is a gonality pencil $g^1_d$ on $C$, then $E_{C,Z}$ is non-simple, and hence, there exists the Donagi-Morrison's extension as in Lemma 5.1;
$$0\longrightarrow M\longrightarrow E_{C,Z}\longrightarrow N\otimes\mathcal{I}_{Z^{'}}\longrightarrow 0.$$
\noindent Since $M\in\mathcal{A}(H)$, $M.N\leq d$, and $\Cliff(C)=2$, we have $M.N=4$ or 5. Hence, we have $M.H\neq N.H$. In fact, since $N^2$ and $M^2$ are even, if $M.H=N.H=9$, then we have $M.N=5$, and hence, we have
$$H.(M\otimes N^{\vee})=0\text{ and }(M\otimes N^{\vee})^2=-2.$$
\noindent However, this contradicts to the ampleness of $H$. Moreover, since $3|M.H$ and $3|N.H$, by the ampleness of $H^{'}=\pi^{\ast}\mathcal{O}_{\mathbb{P}^2}(1)$, we have
$$(M.H,N.H)=(6,12)\text{ or }(12,6).$$
\noindent Hence, we have $M.N=4$. 

If $d=4$, we have $\length Z^{'}=0$, and hence, by the same reason as in the proof of Proposition 6.2, we have
$$E_{C,Z}=\pi^{\ast}\mathcal{O}_{\mathbb{P}^2}(1)\oplus\pi^{\ast}\mathcal{O}_{\mathbb{P}^2}(2).$$

If $d=5$, since $\length Z^{'}=1$, we have $(M.H,N.H)=(12,6)$ (otherwise, since $h^0(M\otimes N^{\vee})=0$, we have a contradiction, by Lemma 5.1). Hence, by the same reason as above, we have
$$M\cong\pi^{\ast}\mathcal{O}_{\mathbb{P}^2}(2)\text{ and }N\cong\pi^{\ast}\mathcal{O}_{\mathbb{P}^2}(1).$$

\smallskip

\smallskip

\noindent {\bf{Acknowledgements}}. The author would like to thank Prof. Hara who has motivated me to consider the splitting problem for ACM bundles in a new light. The author is partially supported by Grant-in-Aid for Scientific Research (25400039), Japan Society for the Promotion Science.

$\;$

\noindent Kenta Watanabe \thanks{Department of Mathematical Sciences, 

\noindent Osaka University, 1-1 Machikaneyama-chou Toyonaka Osaka 560-0043 Japan, 

\noindent {\it E-mail address:goo314kenta@mail.goo.ne.jp}, 

\noindent Telephone numbers: 090-9777-1974}

\end{document}